\font\we=cmb10 at 14.4truept
\font\li=cmb10 at 12truept
~\vskip 1.0cm
\centerline {\we Refined Brill-Noether Locus}
\centerline {\bf and}
\centerline {\we Non-Abelian Zeta
Functions  for Elliptic  Curves}
\vskip 1.0cm
\centerline {\bf Lin WENG}
\centerline {\bf Graduate School of Mathematics, Nagoya University, Japan}
\vskip 0.45cm
In this paper,  new local and global non-abelian zeta functions for
elliptic curves are defined using moduli spaces of semi-stable bundles.
To understand them, we also introduce and study certain refined
Brill-Noether locus in the moduli spaces. Examples of these new zeta
functions and a justification of using only semi-stable bundles are given
too. We end this paper with an appendix on the so-called Weierstrass
groups for general curves, which is motivated by a construction of
Euler systems using torsion points of elliptic curves.
\vskip 0.30cm
\centerline {\li 1. Refined Brill-Noether Locus}
\vskip 0.30cm
\noindent
{\bf 1.1. Moduli Space of Semi-Stable Bundles}

\noindent
{\it 1.1.1 Indecomposable Bundles}

Let $E$ be an elliptic curve defined over $\overline{{\bf F}_q}$,
an algebraic closure of the finite field ${\bf F}_p$ with $q$-elements.

Recall that a vector bundle $V$ on $E$ is called indecomposable if $V$ is
not the direct sum of two proper subbundles, and that every
vector bundle on $E$ may be written as a direct sum of indecomposable
bundles, where the summands and their multiplicities are uniquely
determined up to isomorphism. Thus to understand vector bundles, it
suffices to study the indecomposable ones. To this end, we have the
following result of Atiyah [At]. In the sequel, for simplicity, we always
assume that  the characteristic of ${\bf F}_q$ is strictly bigger than the
rank of $V$.

\noindent
{\bf Theorem.} (Atiyah) {\it (a) For any $r\geq 1$, there is a unique
indecomposable vector bundle $I_r$ of rank $r$ over $E$, all of whose
Jordan-H\"older constituents are isomorphic to ${\cal O}_E$. Moreover,
the bundle $I_r$ has a canonical filtration $$\{0\}\subset
F^1\subset\dots\subset F^r=I_r$$ with $F^i=I_i$ and $F^{i+1}/F^i={\cal
O}_E$;

\noindent
(b) For any $r\geq 1$ and any integer $a$, relative prime to $r$ and each
line bundle $\lambda$ over $E$ of degree $a$, there exists up to
isomorphism a unique indecomposable bundle $W_r(a;\lambda)$ over $E$ of
rank $r$ with $\lambda$ the determinant;

\noindent
(c) The bundle $I_r(W_{r'}(a;\lambda))=I_r\otimes W_{r'}(a;\lambda)$is
indecomposable and every indecomposable bundle is isomorphic to
$I_r(W_{r'}(a;\lambda))$ for a suitable choice of $r,r',\lambda$. Every
bundle $V$ over $E$ is a direct sum of vector bundles of the form
$I_{r_i}(W_{r_i'}(a_i;\lambda_i))$, for  suitable choices of
$r_i,r_i'$ and $\lambda_i$. Moreover, the triples
$(r_i,r_i',\lambda_i)$ are uniquely specified up to permutation by the
isomorphism type of $V$.}

Here note in particular that  $W_r(0,\lambda)\simeq \lambda$, and
that indeed $I_r(W_{r'}(a;\lambda))$ is the unique indecomposable bundle
of rank $rr'$ such that all of whose successive quotients in the
Jordan-H\"older filtration are isomorphic to $W_{r'}(a;\lambda)$.

\noindent
{\it 1.1.2. Semi-Stable Bundles}
 
As above, let $V$ be a vector bundle over $E$. Define its
slop $\mu(V)$ by $\mu(V):={\rm deg}(V)/{\rm rank}(V)$.
Then, following Mumford ([Mu]), $V$ is called stable (resp.
semi-stable), if for any proper subbundle $W$ of $V$,
$\mu(W)<\mu(V)$ (resp. $\mu(W)\leq \mu(V)$). For example,
$I_r(W_{r'}(a;\lambda))$ is semi-stable with
$\mu(I_r(W_{r'}(a;\lambda)))=a/r'$.

\noindent
{\bf Theorem.} (Atiyah) {\it (a) Every bundle $V$ over $E$ is isomorphic
to a direct sum $\oplus_iV_i$ of semi-stable bundles, where
$\mu(V_i)>\mu(V_{i+1})$;

\noindent
(b) Let $V$ be a semi-stable bundle over $E$ with slop $\mu(V)=a/r'$ where
$r'$ is a positive integer and $a$ is an integer relatively prime to
$r'$. Then $V$ is a direct sum of bundles of the form
$I_r(W_{r'}(a;\lambda))$, where $\lambda$ is a line bundle of degree
$a$.}
\vfill\eject
\noindent
{\it 1.1.3. Moduli Space of Semi-Stable Bundles}

Let $V$ be a semi-stable vector bundle, then we may associate it a
Jordan-H\"older filtration, which is far from being unique. However
the associated graded bundle, denoted as ${\rm gr}(V)$, is unique.
Following Seshadri, two semi-stable vector bundles $V$ and
$V'$ are called S-equivalent, denoted by $V\sim_SV'$, if ${\rm
gr}(V)\simeq {\rm gr}(V')$.

Now set $${\cal M}_{E,r}(\lambda)=\{V:{\rm semi}-{\rm stable}, {\rm
det}V=\lambda, {\rm rank}(V)=r\}/\sim_S.$$ Then from the above
classification we have the following well-known

\noindent
{\bf Theorem.} (Atiyah, Mumford-Seshadri) {\it With respect to a fixed
pair $(r,\lambda)$, there exists a natural projective algebraic variety
structure on  ${\cal M}_{E,r}(\lambda)$. Moreover, if $\lambda\in {\rm
Pic}^0(E)$, then
${\cal M}_{E,r}(\lambda)$ is simply the projective space ${\bf
P}^{r-1}_{\overline {\bf F}_q}$.}
\vskip 0.45cm
\noindent
{\bf 1.2. Refined Brill-Noether Locus}

\noindent
{\it 1.2.1. Rational Points}

Now  let $E$ be an elliptic curve defined over a
finite field ${\bf F}_q$. Then over $\overline E=E\times_{{\bf
F}_q}\overline{{\bf F}_q}$, from 1.1.3, we have the
moduli spaces  ${\cal M}_{{\bar E},r}(\lambda)$
(resp. ${\cal M}_{{\bar E},r}(d)$) of semi-stable bundles of
rank $r$ with determinant $\lambda$ (resp. degree $d$) over $\bar E$. As
algebraic varieties, we may consider ${\bf F}_q$-rational points of these
moduli spaces. Clearly, by  definition, these
rational points of moduli spaces correspond exactly to these classes
of semi-stable bundles which themselves are defined over
${\bf F}_q$. (In the case for ${\cal M}_{{\bar E},r}(\lambda)$,
$\lambda$ is assumed to be rational over ${\bf F}_q$.) Thus for simplicity,
we simply write ${\cal M}_{E,r}(\lambda)$ or ${\cal M}_{E,r}(d)$ for
the corresponding  subsets of ${\bf F}_q$-rational points. For
example, we then simply write
${\rm Pic}^0(E)$ for ${\rm Pic}^0(E)({\bf F}_q)$. And by an abuse of
notation, we often call these subsets the moduli spaces of semi-stable
bundles too.

\noindent
{\it 1.2.2. Standard Brill-Noether Locus}

Note that if $V$ is semi-stable with strictly positive degree $d$, then
$h^0(E,V)=d$. Hence the standard Brill-Noether locus is either the whole
space or empty. In this way, we are lead to study the case when $d=0$.

For this, recall that for
$\lambda\in {\rm Pic}^0(E)$,
$${\cal M}_{E,r}(\lambda)=\{V:{\rm semi}-{\rm stable}, {\rm
rank}(V)=r,{\rm det}(V)=\lambda\}/\sim_S$$  is identified with
$$\{V=\oplus_{i=1}^rL_i:\otimes_i L_i=\lambda,L_i\in {\rm Pic}^0(E),
i=1,\dots,r\}/\sim_{\rm iso}\simeq {\bf P}^{r-1}$$ where $/\sim_{\rm iso}$
means modulo isomorphisms.

Now introduce the standard Brill-Noether locus
$$W_{E,r}^a(\lambda):=\{[V]\in {\cal M}_{E,r}(\lambda):h^0(E,{\rm
gr}(V))\geq a\}$$  and its \lq stratification' by
$$W_{E,r}^a(\lambda)^0:=\{[V]\in W_{E,r}(\lambda):h^0(E,{\rm
gr}(V))=a\}=W_{E,r}^a(\lambda)\backslash \cup_{b\geq
a+1}W_{E,r}^b(\lambda).$$
One checks easily that $W_{E,r}^a(\lambda)\simeq {\bf P}^{(r-a)-1}$.
Thus in particular, we have the following

\noindent
{\bf Lemma.} {\it With the same notation as above,
$$W_{E,r+1}^{a+1}(\lambda)\simeq
W_{E,r}^a(\lambda),\qquad{\rm
and}\qquad W_{E,r+1}^{a+1}(\lambda)^0\simeq W_{E,r}^a(\lambda)^0.$$}

\noindent
{\it 1.2.3. Refined Brill-Noether Locus}

The Brill-Noether theory is based on the consideration of $h^0$. But in
the case for elliptic curves, for arithmetic consideration, such a theory
is not fine enough: not only $h^0$ plays important role, the automorphism
groups are important as well. Based on this, we introduce,
for a fixed
$(k+1)$-tuple non-negative integers
$(a_0;a_1,\dots,a_k)$, the subvariety of $W_{E,r}^{a_0}$ by
setting
$$W_{E,r}^{a_0;a_1,\dots,a_k}(\lambda):=\{[V]\in W_{E,r}^{a_0}(\lambda):
{\rm gr}(V)={\cal
O}_E^{(a_0)}\oplus\oplus_{i=1}^kL_i^{(a_i)},\ \otimes_iL_i^{\otimes
a_i}=\lambda,\ L_i\in {\rm Pic}^0(E),\ i=1,\dots,k\}.$$ Moreover, we
define the associated \lq stratification' by setting
$$W_{E,r}^{a_0;a_1,\dots,a_k}(\lambda)^0:=\{[V]\in
W_{E,r}^{a_0;a_1,\dots,a_k}(\lambda),\ \#\{{\cal
O}_E,L_1,\dots,L_k\}=k+1\}.$$
 From definition, we easily have the following

\noindent
{\bf Lemma.} {\it With the same notation as above,
$$W_{E,r+1}^{a_0+1;a_1,\dots,a_k}(\lambda)\simeq
W_{E,r}^{a_0;a_1,\dots,a_k}(\lambda)$$ and
$$W_{E,r+1}^{a_0+1;a_1,\dots,a_k}(\lambda)^0\simeq
W_{E,r}^{a_0;a_1,\dots,a_k}(\lambda)^0.$$ Moreover,
$${\cal
M}_{E,r}(\lambda)=\cup_{a_0;a_1,\dots,a_k}
W_{E,r}^{a_0;a_1,\dots,a_k}(\lambda)^0,$$ where the union is a disjoint
one.}

In fact, the structures of $W_{E,r}^{a_0;a_1,\dots,a_k}(\lambda)$ can be
given explicitly: They are products of (copies of) projective bundles
over $E$ and (copies of) projective spaces.

\noindent
{\bf Proposition.} {\it With the same notation as above, regroup
$(a_0;a_1,\dots,a_k)$ as $(a_0;b_1^{(s_1)},\dots,b_l^{(s_l)})$ with the
condition that
$b_1>b_2>\dots>b_l$ and $s_1,s_2,\dots,s_l\in {\bf Z}_{>0}$, then

\noindent
(1) if $b_l=1$, $$W_{E,r}^{a_0;a_1,\dots,a_k}(\lambda)\simeq
\prod_{i=1}^{l-1}{\bf P}_E^{s_i-1}\times {\bf P}^{s_l};$$

\noindent
(2) if $b_l>1$,
$$W_{E,r}^{a_0;a_1,\dots,a_k}(\lambda)\simeq
\prod_{i=1}^{l}{\bf P}_E^{s_i}.$$}

\noindent
{\it Proof.} This is because we have the following two facts about the
quitient of t products of elliptic curves:

\noindent
(1) The quotient space $E^{(n)}/S_n$ is isomorphic to the ${\bf
P}^{n-1}$-bundle over $E$; and

\noindent
(2) The quotient of $E^{(n-1)}/S_n$ is isomorphic to ${\bf
P}^{(n-1)}$. Here we embed $E^{(n-1)}$ as a subspace of $E^{(n)}$ under
the map:
$$(x_1,\dots,x_n)\mapsto (x_1,\dots,x_{n-1},x_n)$$ with
$x_n=\lambda-(x_1+x_2+\dots+x_{n-1}).$

We end this subsection with the following intersection theoretical
discussion. Fro simplicity, let $\lambda={\cal O}_E$. Then from above, we
have the refined Brill-Noether loci $W_{E,r}^{a_0;a_1,\dots,a_k}({\cal
O}_E)$ which are isomorphic to products of (copies of) projective bundles
over $E$ and (copies of) projective spaces. Thus it would be very
interesting to see the intersections of these special subvarieties in
${\cal M}_{E,r}({\cal O}_E)={\bf P}^{r-1}$. For this purpose, define the
so-called Brill-Noether tautological ring ${\bf BN}_{E,r}({\cal O}_E)$ to
be the subring  generated by all
the associated refined Brill-Noether loci.  For examples,

\noindent
(1) If $r=2$, then this
ring contains two elements: 1-dimensional one $W_{E,2}^{2;0}({\cal
O}_E)=\{[{\cal O}_E\oplus {\cal O}_E]\}$ and the whole ${\bf P}^1$;

\noindent
(2) If $r=3$, then (generators of) this ring contains five elements: 2 of
0-dimensional objects: $W_{E,3}^{3;0}({\cal O}_E)=\{[{\cal O}_E^{(3)}]\}$
and
$W_{E,3}^{1;2}({\cal O}_E)=\{[{\cal O}_E\oplus T_2^{(2)}]:T_2\in E_2\}$
containing 4 elements; 2 of 1-dimensional objects:
$W_{E,3}^{1;1,1}=\{[{\cal O}_E\oplus L\oplus L^{-1}]:L\in {\rm
Pic}^0(E)\}\simeq {\bf P}^1$, a degree 2 projective line contained in
${\bf P}^2={\cal M}_{E,3}({\cal O}_E);$ and
$W_{E,3}^{0;2,1}=\{[L^{(2)}\oplus L^{-2}]:L\in {\rm Pic}^0(E)\}$ a degree
3 curve which is isomorphic to
$E$; and finally the whole space. Moreover, the intersection of
$W_{E,3}^{1;1,1}={\bf P}^1$ and $W_{E,3}^{0;2,1}=E$ are supported on
0-dimensional locus $W_{E,3}^{1;1,1}$, with the multiplicity 3 on the
single point locus $W_{E,3}^{3;0}({\cal O}_E)$ and 1 on the completement
of the points in $W_{E,3}^{1;1,1}$.
\eject
\vskip 0.45cm
\centerline {\li 2. Invariants $\alpha,\beta$ and $\gamma$}
\vskip 0.30cm
\noindent
{\bf 2.1. Measure Refined Brill-Noether Locus Arithmetically}

\noindent
{\it 2.1.1. Invariant $\alpha$}

In the rest of this section, we use the same notation as in 1.2.

To measure the Brill-Noether locus, we introduce the following arithmetic
invariant $\alpha_{E,r}(\lambda)$ by setting
$$\alpha_{E,r}(\lambda):=\sum_{V\in [V]\in {\cal
M}_{E,r}(\lambda)}{{q^{h^0(E,V)}}\over {\#{\rm Aut}(V)}}.$$
Also set
$$\alpha^{a_0+1;a_1,\dots,a_k}_{E,r}(\lambda):=\sum_{V\in [V]\in
W^{a_0+1;a_1,\dots,a_k}_{E,r}(\lambda)^0}{{q^{h^0(E,V)}}\over {\#{\rm
Aut}(V)}}.$$

Before going further, we remark that above, we write $V\in [V]$ in the
summation. This is because in each
$S$-equivalence class
$[V]$, there are usually more than one vector bundles $V$. For example,
$[{\cal O}_E^{(4)}]$ consists of ${\cal O}_E^{(4)}$,  ${\cal
O}_E^{(2)}\oplus I_2$, $I_2\oplus I_2$, ${\cal O}_E\oplus I_3$, and $I_4$
by the result of Atiyah cited in 1.1.

Thus, by Lemma 1.2.3, we have the following

\noindent
{\bf Lemma.} {\it With the same notation as above,
$$\alpha_{E,r}(\lambda)=\sum_{(a_0;a_1,\dots,a_k);k}\alpha^{a_0;a_1,\dots,a_ 
k}_{E,r}(\lambda).$$}

We end this section with the following

\noindent
{\bf Conjecture.} {\it For all $\lambda\in {\rm Pic}^0(E)$,
$$\alpha_{E,r}(\lambda)=\alpha_{E,r}({\cal O}_E).$$}

\noindent
{\it 2.1.2. Invariants $\beta$ and $\gamma$}

Due to the importance of automorphism groups, following
Harder-Narasimhan, and Desale-Ramanan, we introduce the following
$\beta$-series invariants $\beta_{E,r}(d)$, $\beta_{E,r}(\lambda)$
and
$\beta_{E,r}^{a_0;a_1,\dots,a_k}(\lambda)$
by setting
$$\beta_{E,r}(d):=\sum_{V\in [V]\in {\cal
M}_{E,r}(d)}{{1}\over {\#{\rm Aut}(V)}},$$

$$\beta_{E,r}(\lambda):=\sum_{V\in [V]\in {\cal
M}_{E,r}(\lambda)}{{1}\over {\#{\rm Aut}(V)}},$$
and
$$\beta^{a_0;a_1,\dots,a_k}_{E,r}(\lambda):=\sum_{V\in [V]\in
W^{a_0;a_1,\dots,a_k}_{E,r}(\lambda)^0}{{1}\over {\#{\rm
Aut}(V)}}.$$
Corresponding to the Conjecture 2.1.1 for $\alpha$, for
$\beta$, we have the following deep

\noindent
{\bf Theorem.} ([HN] \& [DR]) {\it For all
$\lambda,\lambda'\in {\rm Pic}^d(E)$,
$$\beta_{E,r}(\lambda)=\beta_{E,r}(\lambda').$$
Moreover,
$$N_1\cdot \beta_{E,r}(\lambda)={{N_1}\over
{q-1}}\cdot\prod_{i=2}^r\zeta_E(i)-\sum_{\Sigma_1^kr_i=r,\Sigma_i
d_i=d,{{d_1}\over {r_1}}>\dots>{{d_k}\over{r_k}}, k\geq
2}\prod_i\beta_{E,r_i}(d_i){1\over{q^{\Sigma_{i<j}(r_jd_i-r_id_j)}}}.$$
Here $N_1$ denotes $\#E(:=\#E({\bf F}_q))$ and $\zeta_E(s)$ denotes the
Artin zeta function for elliptic curve $E/{\bf F}_q$.}

Thus, we are lead to introduce the $\gamma$-series invariants
$\gamma_{E,r}(\lambda)$ and
$\gamma_{E,r}^{a_0+1;a_1,\dots,a_k}(\lambda)$ by setting
$$\gamma:=\alpha-\beta.$$ That is to say,
$$\gamma_{E,r}(\lambda):=\sum_{V\in [V]\in {\cal
M}_{E,r}(\lambda)}{{q^{h^0(E,V)}-1}\over {\#{\rm Aut}(V)}},$$
and
$$\gamma^{a_0;a_1,\dots,a_k}_{E,r}(\lambda):=\sum_{V\in [V]\in
W^{a_0;a_1,\dots,a_k}_{E,r}(\lambda)^0}{{q^{h^0(E,V)}-1}\over {\#{\rm
Aut}(V)}}.$$

Clearly, by the above Theorem, the Conjecture 2.1.2 for $\alpha$ is
equivalent to the following

\noindent
{\bf Conjecture.} {\it For all $\lambda\in {\rm Pic}^0(E)$,
$$\gamma_{E,r}(\lambda)=\gamma_{E,r}({\cal O}_E).$$}

The advantage of this Conjecture is that now the support of the
summation is over $W_{E,r}^1(\lambda)$, a codimension 1 projective
subspace.

Similarly, we have the following

\noindent
{\bf Lemma.} {\it With the same notation as above, for $\lambda\in {\rm
Pic}^0(E)$,
$$\beta_{E,r}(\lambda)=\sum_{a_0;a_1,\dots,a_k}\beta^{a_0;a_1,\dots,a_k}_{E, 
r}(\lambda)$$
and
$$\gamma_{E,r}(\lambda)=\sum_{a_0;a_1,\dots,a_k}\gamma^{a_0;a_1,\dots,a_k}_{ 
E,r}(\lambda).$$}

\noindent
{\bf 2.2. Relations Between $\beta$ and $\gamma$}

\noindent
{\it 2.2.1. Bundles with Trivial Graded Bundles}

After certain painful calculations, by an accident, we are lead to the
following

\noindent
{\bf Conjecture.} {\it For all $\lambda\in {\rm Pic}^0(E)$,
$$\gamma^{a_0+1;a_1,\dots,a_k}_{E,r+1}(\lambda)=
\beta^{a_0;a_1,\dots,a_k}_{E,r}(\lambda).$$}
 
Note that by Lemma 1.2.3, $$W_{E,r}^{a_0;a_1,\dots,a_k}(\lambda)^0\simeq
W_{E,r+1}^{a_0+1;a_1,\dots,a_k}(\lambda)^0.$$ Hence, from the definition,
to verify the latest Conjecture, it suffices to show that for any fixed
$(L_1,\dots,L_k)\in \big({\rm Pic}^0(E)\big)^{(k)}$ such that $\#\{{\cal
O}_E,L_1,\dots,L_k\}=k+1$, we have
$$\sum_{V:{\rm gr}(V)={\cal
O}_E^{(a_0)}\oplus\oplus_{i=1}^kL_i^{(a_i)}}{1\over {\#{\rm Aut}(V)}}
=\sum_{V:{\rm gr}(V)={\cal
O}_E^{(a_0+1)}\oplus\oplus_{i=1}^kL_i^{(a_i)}}{{q^{h^0(E,V)}-1}\over
{\#{\rm Aut}(V)}}.$$
Thus by looking at the structure of the automorphism groups carefully
from the  condition that  ${\cal
O}_E,L_1,\dots,L_k$ are all different, we see that this latest conjecture
is equivalent to the following
\vskip 0.30cm
\noindent
{\bf Main Conjecture} (in Algebraic Theory of Non-Abelian Zeta Functions
for Elliptic Curves.) {\it For any $r\in {\bf Z}_{>0}$,
$$\sum_{V:{\rm gr}(V)={\cal O}_E^{(r)}}{1\over {\#{\rm Aut}(V)}}
=\sum_{W:{\rm gr}(W)={\cal O}_E^{(r+1)}}{{q^{h^0(E,W)}-1}\over {\#{\rm
Aut}(W)}}.$$}

Clearly, the advantage of the Main Conjecture is
that only trivial bundle and its extensions are involved.

Before going further, to convince the reader, we check the example with
$r=3$. In this case, there are three possibilities for $V$. That is,
${\cal O}_E^{(3)}$, ${\cal O}_E\oplus I_2$ and $I_3$. One checks
the cardinal numbers of the corresponding automorphism groups are
$(q^3-1)(q^3-q)(q^3-q^2)$,
$(q-1)^2q^3$ and $(q-1)q^2$ respectively. Similarly,
there are five possibilities for $W$. That is,
${\cal O}_E^{(4)}$, ${\cal O}_E^{(2)}\oplus I_2$,
$I_2^{(2)}$, ${\cal O}_E\oplus I_3$ and $I_4$. The cardinal
numbers of the corresponding automorphism groups are
$(q^4-1)(q^4-q)(q^4-q^2)(q^4-q^3)$,
$(q-1)q^3(q^3-q)(q^3-q^2)$, $q^4(q^2-1)(q^2-q)$, $(q-1)^2q^4$ and
$(q-1)q^3$, and their $h^0$ are given by $4,3,2,2,1$ respectively.
Therefore, the Main Conjecture is equivalent to, in case $r=3$,
the following relation:
$$\eqalign{{1\over
{(q^3-1)(q^3-q)(q^3-q^2)}}&+{1\over{(q-1)^2q^3}}+{1\over{(q-1)q^2}}\cr
=&{{q^4-1}\over
{(q^4-1)(q^4-q)(q^4-q^2)(q^4-q^3)}}+{{q^3-1}\over
{(q-1)q^3(q^3-q)(q^3-q^2)}}\cr
&\qquad+{{q^2-1}\over
{q^4(q^2-1)(q^2-q)}}+{{q^2-1}\over{(q-1)^2q^4}}+{{q-1}\over{(q-1)q^3}}.\cr}$$
We leave the routine check of this latest identity to the reader, who may
certainly be amused if doing correctly.

\noindent
{\it 2.2.2. Main Conjecture Implies All Conjectures}

It suffices to imply Conjecture 2.1.2 from the Main Conjecture. For
this, by Lemma 2.1.2, we need to show that
$$\sum_{a_0;a_1,\dots,a_k}\gamma^{a_0;a_1,\dots,a_k}_{E,r}(\lambda)=
\sum_{a_0;a_1,\dots,a_k}\gamma^{a_0;a_1,\dots,a_k}_{E,r}({\cal O}_E).$$
Now by the Main Conjecture, which is equivalent to Conjecture 2.2.1, the
left hand side becomes
$$\sum_{a_0;a_1,\dots,a_k}\beta^{a_0-1;a_1,\dots,a_k}_{E,r-1}(\lambda)$$
while the right hand side is simply
$$\sum_{a_0;a_1,\dots,a_k}\beta^{a_0-1;a_1,\dots,a_k}_{E,r-1}({\cal
O}_E).$$ Note also that $\gamma_{E,r}^{0;a_1,\dots,a_k}(\lambda)=0$.
Therefore the left hand side is simply $\beta_{E,r-1}(\lambda)$ while the
right hand side becomes $\beta_{E,r-1}({\cal O}_E)$. Thus by Theorem 2.1.2,
we complete the proof of the following

\noindent
{\bf Theorem.} {\it Assume that for any $1\leq r'\leq r$,
$$\sum_{V:{\rm gr}(V)={\cal O}_E^{(r')}}{1\over {\#{\rm Aut}(V)}}
=\sum_{W:{\rm gr}(W)={\cal O}_E^{(r'+1)}}{{q^{h^0(E,W)}-1}\over {\#{\rm
Aut}(W)}}.$$ Then $$\alpha_{E,r}(\lambda)=\beta_{E,r-1}({\cal O}_E)$$
for all $\lambda\in {\rm Pic}^0(E)$. In particular,
$$\alpha_{E,r}(0)=N_1\cdot \alpha_{E,r}({\cal O}_E),\qquad
\beta_{E,r}(0)=N_1\cdot \beta_{E,r}({\cal O}_E),\qquad{\rm and}\quad
\gamma_{E,r}(0)=N_1\cdot \gamma_{E,r}({\cal O}_E).$$}

Thus, theoretically,  with the relation $\alpha-\beta=\gamma$, we then can
determine all invariants $\alpha_{E,r},\beta_{E,r}$ and $\gamma_{E,r}$.
(We reminder the reader that in practice the precise formula in Theorem
2.1.2 is hardly useful as there are too many infinite summations
involved.)
\eject
\vskip 0.45cm
\centerline {\li 3. New Non-Abelian Zeta Functions for Elliptic Curves}
\vskip 0.30cm
\noindent
{\bf 3.1. Non-Abelian Local Zeta Functions}

\noindent
{\it 3.1.1. Definition}

Let $E$ be an elliptic curve defined over ${\bf F}_q$, the finite field
with $q$ elements. Then we have the associated (${\bf F}_q$-rational
points of) moduli spaces ${\cal M}_{E,r}(d)$. By definition, the {\it rank
$r$ non-abelian zeta function} $\zeta_{E,r,{\bf F}_q}(s)$ {\it of} $E$ is
defined by setting
$$\zeta_{E,r,{\bf F}_q}(s):=\sum_{V\in [V]\in {\cal M}_{E,r}(d),d\geq
0}{{q^{h^0(E,V)}-1}\over {\#{\rm Aut}(V)}}q^{-s\cdot d(V)},\qquad {\rm
Re}(s)>1.$$ Here $d(V)$ denotes the degree of $V$.

\noindent
{\it Remark.} We call the above infinite sum the rank $r$ non-abelian zeta
function because when $r=1$ the above summation $\zeta_{E,1}(s)$
coincides with the classical Artin zeta function for the elliptic curve
$E$ over ${\bf F}_q$. In other words,  the classical Artin zeta
function $\zeta_E(s)$ may be written as $$\zeta_E(s)=\sum_{L\in {\rm
Pic}^d(E),d\geq 0}{{q^{h^0(E,L)}-1}\over {\#{\rm Aut}(L)}}q^{-s\cdot
d(L)},\qquad {\rm Re}(s)>1,$$ where $d(L)$ denotes the degree of $L$.

\noindent
{\it 3.1.2. Basic Properties}

With the above definition, by a direct yet long calculation, we have, by
setting $t=q^{-s}$ and $Z_{E,r,{\bf F}_q}(t):=\zeta_{E,r,{\bf F}_q}(s)$,
the following

\noindent
{\bf Fundamental Identity.} {\it Let $E$ be an elliptic curve defined
over ${\bf F}_q$. Then for any $r\in {\bf Z}_{>0}$,
$$Z_{E,r,{\bf F}_q}(t):=\gamma_{E,r}(0)\cdot{{P_{E,r,{\bf F}_q}(t)}\over
{(1-t^r)(1-q^rt^r)}}.$$ Here
$$\eqalign{.&P_{E,r,{\bf
F}_q}(t)\cr
=&1+\sum_{i=1}^{r-1}(q^i-1){{\beta_{E,r}(i)}\over
{\gamma_{E,r}(0)}}
\cdot t^i+\big((q^r-1){{\beta(0)}\over{\gamma_{E,r}(0)}}-(q^r+1)\big)
\cdot t^r
+\sum_{i=1}^{r-1}(q^i-1){{\beta_{E,r}(i)}\over{\gamma_{E,r}(0)}}
\cdot q^{r-i}t^{2r-i}+q^rt^{2r}.\cr}$$}

\noindent
{\it Remark}. We in this paper choose  not to write down the detailed
elementary calculations, despite the fact that some of them are very long
and sometimes a bit complicated.

As a direct consequence of this Fundamental Identity, we have the following

\noindent
{\bf Theorem.} {\it Let $E$ be an elliptic curve defined over ${\bf F}_q$, 
the finite field
with $q$ elements. Then the associated rank $r$ non-abelian zeta function
$\zeta_{E,r,{\bf F}_q}(s)$ satisfies the following basic properties:

\noindent
(1) ({\bf Rationality}) $Z_{E,r,{\bf F}_q}(t)$ may be written as the
quotient of two polynomials;

\noindent
(2) ({\bf Functional Equation}) $\zeta_{E,r,{\bf F}_q}(s)=\zeta_{E,r,{\bf
F}_q}(1-s).$}

 From here we may conclude that after suitable arrangement, the product of
two reciprocal roots of $P_{E,r,{\bf F}_q}(t)$, a degree $2r$
polynomial, are always equal to
$q$, the cardinal number of the base field. Moreover, we know that up to
the term $\gamma_{E,r}(0)$, from Theorem 2.1.2, the coefficients of these
local non-abelian zeta functions can be computed. Thus if  the Main
Conjecture is assumed, then  all terms of our
non-abelian zeta functions can be given precisely.

Moreover, note that the moduli spaces ${\cal M}_{E,r}(d)$ are indeed
projective bundles over $E$. Thus if the main conjecture holds, in the
study of non-abelian zeta functions for elliptic curves, we are lead to
study only the  vertical direction, i.e., the Brill-Noether loci
appeared a single fiber (and hence all fibers). So, the horizontal
direction, i.e., the Picard group, plays no role. This is why we call the
elegant conjectural relation
$$\sum_{V:{\rm gr}(V)={\cal O}_E^{(r)}}{1\over {\#{\rm Aut}(V)}}
=\sum_{W:{\rm gr}(W)={\cal O}_E^{(r+1)}}{{q^{h^0(E,W)}-1}\over {\#{\rm
Aut}(W)}}$$  the Main Conjecture in algebraic theory of non-abelian zeta
functions for elliptic curves.

\vskip 0.30cm
\noindent
{\bf 3.2. Non-Abelian Global Zeta Functions}

\noindent
{\it 3.2.1. Definition.}

Now let ${\bf E}$ be an elliptic curve defined over
a number field $F$. For a fixed positive integer $r$, set $S_{\rm bad}$ to
be the union of all infinite places (of $F$), all finite places
where ${\bf E}$ have bad reductions, or where
the characteristics of the residue fields are less than $r$. By
definition, a (finite) places $v$ of $F$ is good, if $v$ is not in $S_{\rm
bad}$.

Thus, in particular, for good places $v$, by taking reduction of ${\bf E}$
at $v$, we have the associated regular elliptic curve ${\bf E}_v$ defined
over
$F(v)\simeq {\bf F}_{q_v}$, the residue field of $F$ at $v$, where $q_v$
denotes cardinal number of $F(v)$. Then by applying the construction of
3.1, we get the rank
$r$ local non-abelian zeta functions $\zeta_{{\bf E}_v,r,{\bf
F}_{q_v}}(s)$. In particular, we further obtain, by the rationality,
the corresponding polynomials $P_{{\bf E}_v,r, {\bf F}_q}(t)$ of degree
$2r$ (with 1 as  constant terms).

\noindent
{\bf Definition.} {\it Let ${\bf E}$ be an elliptic curve defined over a
number field $F$. Then for any positive integer $r$, define its associated
rank $r$ global non-abelian zeta function $\zeta_{{\bf E},r,F}(s)$ by
setting
$$\zeta_{{\bf E},r,F}(s):=\prod_{v:{\rm good}}P_{{\bf
E}_v,r,{\bf F}_{q_v}}(q_v^{-s})^{-1}.$$ Here $q_v$ denotes the cardinal
number of the residue field of $F$ at $v$.}
 
Clearly, if $r=1$, this then recovers the famous
Hasse-Weil zeta function for elliptic curves, for which we have the
celebrated BSD conjecture. Thus it seems to be quite natural for us to
call the above Euler product a global non-abelian zeta function for
elliptic curve.

Surely the biggest problem we are now facing in this algebraic part of
our non-abelian zeta function is to give the precise region over which
the Euler product in the definition converges. For this we have to use
our refined Brill-Noether theory discussed in Section 1.

\noindent
{\it 3.2.2. Estimations for $\beta$ and $\gamma$}

We now want to give  estimations for invariants $\beta$ and $\gamma$. So
let $E$ be an elliptic curve defined over a finite field
${\bf F}_q$ as before.

First we study $\beta_{E,r}$. For this, following Harder-Narasimhan,
we interpret the Tamagawa number of ${\rm SL}(n)$ is 1  as follows:

\noindent
{\bf Proposition.} (=[DR Proposition 1.1])
{\it Let $\zeta_E(s)$ be the Artin zeta function of $E$. Then for any
fixed $\lambda\in {\rm Pic}^d(E)$,
$$\sum_{V:{\rm rank}(V)=r,{\rm det}(V)=\lambda}{1\over {\#{\rm Aut}(V)}}=
{1\over {q-1}}\prod_{k=2}^r\zeta_E(k).$$}

Therefore, for a fixed $\lambda\in {\rm Pic}^d(E)$, $$\sum_{V:{\rm
rank}(V)=r,{\rm det}(V)=\lambda}{1\over {\#{\rm Aut}(V)}}=O(q^{-1}).$$
This implies in particular that
$$\beta_{E,r}(\lambda)=O(q^{-1}).$$ Thus, by Hasse's result ([Ha]) on
$N_1=\#{\rm Pic}^d(E)$, we know $N_1=O(q)$. Therefore we complete the
following

\noindent
{\bf Proposition I.} {\it With the same notation as above,
$\beta_{E,r}(d)=O(1)$, as $q\to\infty$.}

Next we study $\gamma_{E,r}(0)$. As our final purpose is to give a good
estimation for the coefficients of $P_{E,r,{\bf F}_q}(t)$, hence by the
Fundamental Identity in 3.1.2, it suffices to give a lower bound for
$\gamma_{E,r}(0)$. For this, we consider semi-stable vector bundles $V$
with ${\rm gr}(V)={\cal O}_E\oplus\oplus_{i=1}^{r-1}L_i$ with
$L_i\in {\rm Pic}^0(E)$ and $\#\{{\cal O}_E,L_1,\dots,L_r\}=r$. Clearly
then $V={\rm gr}(V)$ and there are totally $O(N_1^{r-1})$ or better
$O(q^{r-1})$ of them, as $q\to\infty$ by the above mentioned result of
Hasse. On the other hand, easily, we have $h^0(V)=1$ and $\#{\rm
Aut}(V)=(q-1)^r$. Therefore
$$\sum_{V:{\rm gr}(V)={\cal O}_E\oplus\oplus_{i=1}^{r-1},L_i\in {\rm
Pic}^0(E),\#\{{\cal O}_E,L_1,\dots,L_r\}=r+1}{{q^{h^0(E,V)}-1}\over
{\#{\rm Aut}(V)}}=O(1).$$ This then implies the following

\noindent
{\bf Proposition II.} {\it With the same notation as above,
$\gamma_{E,r}(d)=O(1)$, as $q\to\infty$.}

\noindent
{\it Remark.} In fact if we assume the Main Conjecture in 2.2.1, i.e.,
$$\sum_{V:{\rm gr}(V)={\cal O}_E^{(r)}}{1\over {\#{\rm Aut}(V)}}
=\sum_{W:{\rm gr}(W)={\cal O}_E^{(r+1)}}{{q^{h^0(E,W)}-1}\over {\#{\rm
Aut}(W)}},$$ then by Theorem 2.2.2, $\gamma_{E,r}(0)=\beta_{E,r-1}(0)$.
Hence Proposition II is a direct consequence of Proposition I.

\noindent
{\it 3.2.3. Convergence of Global Non-Abelian Zeta Functions}

As above, let $E$ be an elliptic curve defined over ${\bf F}_q$. Then by
the Fundamental Identity and the Functional Equation for local
non-abelian zeta functions of elliptic curves,
$$P_{E,r,{\bf F}_q}(t)=\prod_{i=1}^r(1+A_it+qt^2)$$ with $A_i\in {\bf
R}.$ Therefore, we have
$$\sum_iA_i=(q-1){{\beta_{E,r}(1)}\over{\gamma_{E,r}(0)}};\qquad
\sum_{i<j}A_iA_j=(q^2-1){{\beta_{E,r}(2)}\over{\gamma_{E,r}(0)}};$$
$$\dots,\qquad
\prod_{i=1}^rA_i=(q^r-1){{\beta_{E,r}(0)}\over{\gamma_{E,r}(0)}}-(q^r+1).$$
So, as $q\to\infty$, by Proposition I and II in 3.2.1,
$$\eqalign{\sum_iA_i=&O(q);\cr
\sum_{i<j}A_iA_j=&O(q^2);\cr
\dots\dots&\dots\dots\dots\cr
\prod_{i=1}^rA_i=&O(q^r).\cr}$$
So, $A_i=O(q), i=1,\dots,r$. As a direct consequence, we then have
the following
\vskip 0.30cm
\noindent
{\bf Theorem.} {\it Let ${\bf E}$ be an elliptic curve defined over a
number field $F$. Then all rank $r$ global non-abelian zeta function
$\zeta_{{\bf E},r,F}(s)$, defined by (infinite) Euler products converge
when ${\rm Re}(s)>2$.}

\noindent
{\it Remark.} When $r=1$, by the above mentioned result of Hasse, the
Hasse-Weil zeta functions, i.e., $\zeta_{{\bf E},1,F}(s)$ in our notation,
converge when
${\rm Re}(s)>{3\over 2}$. While this result is better than ours in the
case when $r=1$,  we are also quite satisfied with our one as our Theorem
is the best possible for general
$r$ at this stage. (See e.g. 4.1.1 below.)

Now, a natural question is whether our rank $r$ global non-abelian
zeta functions have meromorphic extensions and satisfy the functional
equation. We believe that it should be the case. In fact we have the
following

\noindent
{\bf Working Hypothesis.} {\it By introducing also factors for bad
places, the completed rank $r$ non-abelian zeta functions $\xi_{{\bf
E},r,F}(s)$ for elliptic curves have meromorphic continuation to the whole
complex plane and satisfy the functional equation
$$\xi_{{\bf
E},r,F}(s)=\xi_{{\bf
E},r,F}(1+{1\over r}-s).$$}

Clearly, we then would also hope for certain type of such zeta functions,
the inverse Mellin transform would lead to modular forms of
fractional weight $1+{1\over r}$. Unfortunately, we have not yet obtained
any examples to support this speculation. But if it holds, we then have a
systematical way to construct fractional weight modular forms, for which,
except in the case of half integers we know very little.
\vskip 0.45cm
\centerline {\li 4. Examples and Justifications}
\vskip 0.30cm
\noindent
{\bf 4.1. Lower Rank Non-Abelian Zeta Functions}

\noindent
{\it 4.1.1. Rank Two}

Let $E$ be an elliptic curve defined over the finite field ${\bf F}_q$.
If rank $r$ is two, we need then only to calculate
$\beta_{E,2}(0),\beta_{E,2}(1)$ and $\gamma_{E,r}(0)$.

We first consider $\beta_{E,2}(0)$. Then by our discussion on
Brill-Noether locus, it suffices to calculate
$\beta_{E,2}({\cal O}_E)$.
Now $${\cal M}_{E,2}({\cal O}_E)=W_{E,2}^{2;0}({\cal
O}_E)^0\cup W_{E,2}^{0;2}({\cal O}_E)^0\cup W_{E,2}^{0;1,1}({\cal
O}_E)^0.$$ Clearly, $$W_{E,2}^{2;0}({\cal
O}_E)^0=\{[V]: {\rm gr}(V)={\cal O}_E^{(2)}\}$$ consisting of just 1
element;
$$W_{E,2}^{0;2}({\cal O}_E)^0=\{[V]: {\rm gr}(V)=T_2^{(2)},T_2\in
E_2, T_2\not={\cal O}_E\}$$ consisting of 3 elements coming from
non-trivial $T_2\in E_2$, 2-torsion subgroup of $E$; while
$$W_{E,2}^{0;1,1}({\cal O}_E)^0= \{[V]: {\rm gr}(V)=L\oplus L^{-1},L\in
{\rm Pic}^0(E),L\not=L^{-1}\}$$ is simply the complement
of the above 4 points in ${\bf P}^1$.
With this, one checks that
$$\beta_{E,2}(0)=\Big({1\over {(q^2-1)(q^2-q)}}+{1\over
{(q-1)q}}\Big)+3\cdot \Big({1\over {(q^2-1)(q^2-q)}}+{1\over
{(q-1)q}}\Big)+ \big(q+1-(3+1)\big)\cdot {1\over
{(q-1)^2}}
={{q+3}\over {q^2-1}}.$$ And hence
$$\beta_{E,2}(0)=N_1\cdot {{q+3}\over {q^2-1}}.$$

As for $\beta_{E,2}(1)$, it is very simple: Any degree one rank two
semi-stable bundle is stable. Moreover, by the result of Atiyah cited in
1.1,  there is exactly one stable
rank two bundle whose determinant is the fixed line bundle. Thus
$$\beta_{E,2}(1)=N_1\cdot{1\over {q-1}}.$$

Finally, we study $\gamma_{E,2}(0)$. We want to check Conjecture 2.1.2.
Clearly if $\lambda\not={\cal O}_E$, then $\gamma_{E,2}(\lambda)$ is
supported on $$W^{1;1}_{E,2}(\lambda)=\{[V]:{\rm gr}(V)={\cal
O}_E\oplus\lambda\}$$ consisting only one element with $V={\rm gr}(V)={\cal
O}_E\oplus\lambda.$ So
$$\gamma_{E,2}(\lambda)={{q-1}\over {(q-1)^2}}={1\over {q-1}}.$$

On the other hand,
$\gamma_{E,2}({\cal O}_E)$ is
supported on $$W^{2;0}_{E,2}({\cal O}_E)=\{[V]:{\rm gr}(V)={\cal
O}_E^{(2)}\}$$ consisting only one element too.
However, now in the single class $[V]$ with ${\rm gr}(V)={\cal
O}_E^{(2)}$, there are two elements, i.e., ${\cal
O}_E^{(2)}$ and $I_2$. So
$$\gamma_{E,2}({\cal O}_E)={{q^2-1}\over
{(q^2-1)(q^2-q)}}
+{{q-1}\over{(q-1)q}}={1\over {q-1}}=\beta_{E,1}({\cal O}_E).$$
Thus we have checked all conjectures in the case $r=2$.

\noindent
{\bf Proposition.} {\it With the same notation as above, in the case
$r=2$, all conjectures in this paper are confirmed. In particular,
$$Z_{E,2, {\bf F}_q}(t)={{N_1}\over
{q-1}}\cdot{{1+(q-1)t+(2q-4)t^2+(q^2-q)t^3+q^2t^4}\over
{(1-t^2)(1-q^2t^2)}}.$$}

We reminder the reader that in this case, $P_{E,2,{\bf F}_q}(t)$ is
independent of $E$ and with integer coefficients. Thus for global rank
two non-abelian zeta function, we obtain an absolute Euler product, say in
the case when the base field is {\bf Q},
$$E_2(s)=\prod_{p\geq 3,{\rm prime}}{1\over
{1+(p-1)p^{-s}+(2p-4)p^{-2s}+(p^2-p)p^{-3s}+p^2p^{-4s}}},\qquad {\rm
Re}(s)>2.$$
Clearly we expect more from such a beautiful Euler product.

\noindent
{\it 4.1.2. Rank Three}

First, we check whether the Main Conjecture in 2.1.2 holds. That is to
say, we should show $$\sum_{V,{\rm gr}(V)={\cal O}_E^{(2)}}{1\over
{\#{\rm Aut}(V)}}=\sum_{W,{\rm gr}(W)={\cal
O}_E^{(3)}}{{q^{h^0(E,V)}-1}\over {\#{\rm Aut}(V)}}.$$ By the fact that
${\rm Aut}({\cal O}_E\oplus I_2)=(q-1)^2q^3$, the above identity is
equivalent to
$${1\over {(q^2-1)(q^2-q)}}+{1\over {(q-1)q}}={{q^3-1}\over
{(q^3-1)(q^3-q)(q^3-q^2)}}+{{q^2-1}\over {(q-1)^2q^3}}+{{q-1}\over
{(q-1)q^2}},$$ which may be directly checked. This together with the
similar relation for $r=1$  discussed in 4.1.1 leads to
$$\gamma_{E,3}(0)=N_1\cdot \gamma_{E,3}({\cal
O}_E)=N_1\cdot\beta_{E,2}({\cal O}_E)=N_1\cdot{{q+3}\over {q^2-1}}.$$

So we are left to study $\beta_{E,3}(d)$, $d=0,1,2$. Easily, we have
$$\beta_{E,3}(1)=\beta_{E,3}(2)=N_1\cdot {1\over q-1}$$ since here
all semi-stable bundles become stable. Thus we are lead to consider only
$\beta_{E,3}(0)$. So it suffices to give $\beta_{E,3}(\lambda)$ for any
$\lambda\not={\cal O}_E$. (Despite the fact that
$\beta_{E,r}(\lambda)=\beta_{E,r}({\cal O}_E)$ for any $\lambda\in {\rm
Pic}^0(E)$, in practice, the calculation of
$\beta_{E,r}(\lambda)$ with $\lambda\not={\cal O}_E$ is  easier than
that for $\beta_{E,r}({\cal O}_E)$.)

Now $${\cal M}_{E,3}(\lambda)=\Big((W_{E,3}^{2;1}(\lambda)^0)\cup
W_{E,3}^{1;2}(\lambda)^0\cup W_{E,3}^{1;1,1}(\lambda)^0\Big)\cup
W_{E,3}^{0;3}(\lambda)^0
\cup W_{E,3}^{0;2,1}(\lambda)^0\cup W_{E,3}^{0;1,1,1}(\lambda)^0.$$
Moreover, we have

\noindent
(1) $W_{E,3}^{2;1}(\lambda)^0$ consists a single class [V], i.e., the one
with
${\rm gr}(V)={\cal O}_E^2\oplus\lambda$, which contains two vector bundles,
i.e., ${\cal O}_E^2\oplus\lambda$ and $I_2\oplus\lambda$;

\noindent
(2) $W_{E,3}^{2;1}(\lambda)^0\cup
W_{E,3}^{1;2}(\lambda)^0\cup W_{E,3}^{1;1,1}(\lambda)^0\simeq {\bf P}^1$
with
$W^{1;2}(\lambda)^0$ consists of 4 classes $[V]$, i.e., these such that
${\rm gr}(V)={\cal O}_E\oplus \big(\lambda^{1\over 2}\big)^{(2)}$, where
$\lambda^{1\over 2}$ denotes any of the four square roots of $\lambda$.
Clearly then in each class $[V]$, there are also two vector bundles
${\cal O}_E\oplus \big(\lambda^{1\over 2}\big)^{(2)}$ and ${\cal
O}_E\oplus I_2\otimes \lambda^{1\over 2}$;

\noindent
(3) $W_{E,3}^{0;3}(\lambda)^0
\cup W_{E,3}^{0;2,1}(\lambda)^0\cup W_{E,3}^{0;1,1,1}(\lambda)^0={\bf
P}^2\backslash {\bf P}^1$.

\noindent
(3.a) $W_{E,3}^{0;3}(\lambda)^0$ consists of 9 classes $[V]$,
i.e., these $[V]$ with ${\rm gr}(V)=\big(\lambda^{1\over 3}\big)^{(3)}$
where
$\lambda^{1\over 3}$ denotes any of the 9 triple roots of $\lambda$.
Moreover, in each $[V]$, there are three bundles, i.e.,
$\big(\lambda^{1\over 3}\big)^{(3)}$, $\lambda^{1\over
3}\oplus I_2\otimes \lambda^{1\over 3}$ and $I_3\otimes\lambda^{1\over
3}$.

\noindent
(3.b)  $\Big(W_{E,3}^{2;1}(\lambda)^0\cup
W_{E,3}^{1;2}(\lambda)^0\Big)\cup\Big(W_{E,3}^{0;3}(\lambda)^0\cup
W_{E,3}^{0;2,1}(\lambda)^0\Big)$ is isomorphic to $E$. Moreover, each class
$[V]$ in $W^{0;2,1}(\lambda)^0$ consists of two bundles, i.e.,
$L^{(2)}\oplus \lambda\otimes L^{-2}$ and
$I_2\otimes L\oplus \lambda\otimes L^{-2}$ when ${\rm gr}(V)=L^{(2)}\oplus
\lambda\otimes L^{-2}$.

(One checks that in fact the refined Brill-Noether loci ${\bf P}^1$ and
$E$ appeared above are embedded in ${\bf P}^2$ as degree 2 and 3
regular curves. And hence the intersection should be 6: The
intersection points are at
$[V]$ with ${\rm gr}(V)={\cal O}_E^{(2)}\oplus \lambda$ with multiplicity
2, and ${\cal O}_E\oplus \big(\lambda^{1\over 2}\big)^{(2)}$
corresponding to four square roots of $\lambda$ with multiplicity one.
That is to say, the intersection actually are supported on 
$W_{E,3}^{2;1}(\lambda)^0\cup
W_{E,3}^{1;2}(\lambda)^0$. So it would be very interesting in general
to study the intersections of the refined Brill-Noether loci  as well.)

 From this analysis, we conclude that
$$\eqalign{\beta_{E,3}(\lambda)=&\Big({1\over {(q^2-1)(q^2-q)(q-1)}}
+{1\over {(q-1)q(q-1)}}\Big)\cr
&+4\Big({1\over
{(q^2-1)(q^2-q)(q-1)}}+{1\over
{(q-1)q(q-1)}}\Big)+(q-4)\cdot\Big({1\over
{(q-1)^3}}\Big)\cr &+9\Big({1\over {(q^3-1)(q^3-q)(q^3-q^2)}}+{1\over
{(q-1)^2q^3}}+{1\over {(q-1)q^2}}\Big)\cr
&+\big(N_1-(9+4+1)\big)\cdot\Big({1\over
{(q-1)(q^2-1)(q^2-q)}}+{1\over {(q-1)(q-1)q}}\Big)\cr
&+\big(q^2-(N_1-4-1)\big)\cdot{1\over {(q-1)^3}}.\cr}$$
Clearly, it is not of our best interests to write down the associated zeta
function very precisely.
Still,  in this case, all conjectures in this paper have been confirmed.
In particular, we conclude that $r=2$ is the only case
that the global zeta functions are independent of elliptic curves and
that the polynomials are with integral coefficients instead of rational
coefficients in general.

Also we would reminder the reader that here in fact 2- and 3- torsion
points are involved naturally in the calculation.
\vskip 0.30cm
\noindent
{\bf 4.2. Why  Use only Semi-Stable Bundles}

\noindent
{\it 4.2.1. Degree 0}

At the first glance, Theorem 2.1.2 and Proposition 3.2.2 suggest that in
the definition of non-abelian zeta functions we should consider all vector
bundles, just as what happens in the theory of automorphic $L$-functions.
However, we here use an example with $r=2$ to indicate the opposite.

Thus we frist introduce a new  zeta
function $\zeta_{E,r}^{\rm all}(s)$  by
$$\zeta_{E,r}^{\rm all}(s):=\sum_{V: {\rm
rank}(V)=2}{{q^{h^0(V)}-1}\over {\#{\rm Aut}(V)}}\cdot q^{-sd(V)}.$$ Then
by our discussion on the non-abelian zeta functions associated to
semi-stable bundle, we only need to consider the contribution of rank 2
bundles which are not semi-stable.

We start with a discussion on extension of bundles. Assume that $V$ is
not semi-stable of rank 2. Let $L_2$ be the line subbundle of $V$ with
maximal degree, then
$V$ is obtained from the extension of $L_2:=V/L_1$ by $L_1$ $$0\to L_1\to
V\to L_2\to 0.$$ But $V$ is not semi-stable implies that all such
extensions are trivial.  Thus  $V=L_1\oplus L_2$.
For later use, set $d_i$ to be the degree of
$L_i, i=1,2$. Then $d_1+d_2=d$ the degree of $V$, and
$$\#{\rm Aut}(V)=(q-1)^2\cdot q^{h^0(E,L_1\otimes L_2^{-1})}=(q-1)^2\cdot
q^{d_1-d_2}.$$

Next we study the the contribution of degree 0 vector bundles of rank 2
which are not semi-stable. Note that the support of the summation should
have non-vanishing $h^0$. Thus  $V=L_1\oplus L_2$ where
$L_1\in {\rm Pic}^{d_1}(E)$ with $d_1>0$.
So the contributions of these bundles are  given by
$$\eqalign{\zeta_{E,2}^{=0}(s)=&Z_{E,2}^{=0}(t)\cr
=&\sum_{d=1}^\infty\sum_{L_1\in
{\rm Pic}^d(E),L_2\in {\rm
Pic}^{-d}(E)}{{q^{h^0(L_1)}-1}\over{(q-1)^2q^{h^0(L_1\otimes L_2^{\otimes
-1})}}} ={{N_1^2}\over
{(q-1)^2}}\cdot\sum_{d=1}^\infty{{q^d-1}\over{q^{2d}}}\cr
=&{{qN_1^2}\over{(q^2-1)(q-1)^2}}.\cr}$$

\noindent
{\it 4.2.2. Degree $>0$}

Now we consider all degree strictly positive rank 2 vector bundles which
are not semi-stable. From above we see that $V=L_1\oplus L_2$ with
$d_1>d_2$.
Thus for $h^0(E,V)$, there are three  cases:

\noindent
(i) $d_2>0$, clearly then $h^0(E,V)=d$;

\noindent
(ii) $d_2=0$. Here there are two subcases, namely, (a) if $L_2={\cal
O}_E$, then $h^0(E,V)=d_1+1$; (b) If $L_2\not={\cal O}_E$, then
$h^0(E,V)=d_1$;

\noindent
(iii) $d_2<0$. Then $h^0(E,V)=d_1$.

Therefore, all in all the contribution of strictly positive degree rank 2
bundles which are not semi-stable to the zeta function
$\zeta_{E,r}^{\rm all}(s)$ is given by
$$\zeta_{E,2}^{>0}(s)=Z_{E,2}^{>0}(t)=\Big(\sum_{(i)}+\sum_{(ii.a)}+\sum_{(i 
i.b)}
+\sum_{(iii)}\Big){{q^{h^0(V)}-1}\over {\#{\rm Aut}(V)}}t^d$$ where
$\sum_{(*)}$ means the summation is taken for all vector bundles in case
(*).
 
Hence, we have
$$\eqalign{\sum_{(i)}{{q^{h^0(V)}-1}\over {\#{\rm
Aut}(V)}}=&N_1^2\cdot\sum_{d=1}^\infty\sum_{d_1+d_2=d,
d_1>d_2>0}{{q^d-1}\over{(q-1)^2q^{d_1-d_2}}}t^d,\cr
\sum_{(ii.a)}{{q^{h^0(V)}-1}\over
{\#{\rm
Aut}(V)}}=&N_1\cdot\sum_{d=1}^\infty{{q^{d+1}-1}\over{(q-1)^2q^{d}}}
t^{d},\cr
\sum_{(ii.a)}{{q^{h^0(V)}-1}\over {\#{\rm
Aut}(V)}}=&N_1(N_1-1)\cdot\sum_{d=1}^\infty
{{q^{d}-1}\over{(q-1)^2q^{d}}}t^d,\cr
\sum_{(iii)}{{q^{h^0(V)}-1}\over
{\#{\rm Aut}(V)}}=&N_1^2\cdot\sum_{d=1}^\infty\sum_{d_1+d_2=d,
d_1>0>d_2}{{q^{d_1}-1}\over{(q-1)^2q^{d_1-d_2}}}t^d.\cr}$$

By a direct calculation, we find that
$$\eqalign{\sum_{(i)}{{q^{h^0(V)}-1}\over {\#{\rm
Aut}(V)}}=&{{N_1^2t^3}\over {q-1}}\cdot
{{q^2+q+1+q^2t}\over{(1-t^2)(1-q^2t^2)(q-t)}},\cr
\sum_{(ii.a)}{{q^{h^0(V)}-1}\over
{\#{\rm
Aut}(V)}}=&{{N_1t}\over {q-1}}\cdot{{q+1-t}\over{(q-t)(1-t)}},\cr
\sum_{(ii.a)}{{q^{h^0(V)}-1}\over {\#{\rm
Aut}(V)}}=&{{N_1(N_1-1)}\over {q-1}}\cdot {t\over {(q-t)(1-t)}},\cr
\sum_{(iii)}{{q^{h^0(V)}-1}\over
{\#{\rm
Aut}(V)}}=&{{N_1^2t}\over{(q-1)^2(q^2-1)}}\cdot{{q^2+q-1-qt}\over
{(1-t)(q-t)}}.\cr}$$

\noindent
{\it 4.2.3. Degree $<0$}

Finally we consider the contribution of bundles with strictly negative
degree. First we have the following classification according to
$h^0(E,V)$.

\noindent
(i) $d_1>0>d_2$. Then $h^0(E,V)=d_1$;

\noindent
(ii) $d_1=0>d_2$. Here two subcases. (a) $L_1={\cal O}_E$,
then $h^0(E,V)=1$; (b) $L_1\not={\cal O}_E$, then $h^0(V)=0$;

\noindent
(iii) $0>d_1>d_2$. Here $h^0(V)=0$.

Thus note that the support of $h^0(E,V)$ is only on the cases (i) and
(ii.a), we see that similarly as before,
the contribution of strictly positive degree rank 2 bundles
which are not semi-stable to the zeta function is given by
$$\zeta_{E,2}^{<0}(s)=Z_{E,2}^{<0}(t)=\Big(\sum_{(i)}+\sum_{(ii.a)}\Big){{q^ 
{h^0(V)}-1}\over
{\#{\rm Aut}(V)}}t^d.$$

Hence, we have
$$\eqalign{\zeta_{E,2}^{<0}(s)=Z_{E,2}^{<0}(t)
=&
N_1^2\sum_{d=-1}^{-\infty}\sum_{d_1>0>d_2,
d_1+d_2=d}{{q^{d_1}-1}\over{(q-1)^2q^{d_1-d_2}}}t^d+N_1\cdot\sum_{d=-1}^{-\i 
nfty}
{{q-1}\over {(q-1)^2q^{-d}}}t^d\cr
=&{{N_1^2}\over {(q-1)^2}}\cdot {q\over {(qt-1)(q^2-1)}}
+{{N_1}\over{q-1}}\cdot {1\over {qt-1}}.\cr}$$

I hope now the reader is fully convinced that our definition of
non-abelian zeta function by using moduli space of semi-stable bundles is
much better than that of others: Not only our semi-stable zeta functions
have much neat structure, we also have well-behavior geometric and hence
arithmetic spaces ready to use. In a certain sense, we think the picture
of our non-abelian zeta function is quite similar to that the so-called
new forms: Only after  removing these not-semi-stable contributions, we
can see the intrinsic beautiful structures.
\vskip 0.45cm
\centerline {\bf REFERENCES}
\vskip 0.30cm
\item{[A]} E. Artin, Quadratische K\"orper im Gebiete der h\"oheren
Kongruenzen, I,II, {\it Math. Zeit}, {\bf 19} 153-246 (1924) (See also
{\it Collected Papers}, pp. 1-94,  Addison-Wesley 1965)
\vskip 0.30cm
\item{[At]} M. Atiyah, Vector Bundles over an elliptic curve, {\it
Proc. LMS, VII}, 414-452 (1957) (See also {\it Collected Works}, Vol. 1,
pp. 105-143, Oxford Science Publications, 1988)
\vskip 0.30cm
\item{[DR]} U.V. Desale \& S. Ramanan, Poincar\'e polynomials of the variety of
stable bundles, Maeh. Ann {\bf 216}, 233-244 (1975)
\vskip 0.30cm
\item{[HN]} G. Harder \& M.S. Narasimhan, On the cohomology groups of 
moduli spaces
of vector bundles over curves, Math Ann. {\bf 212}, (1975) 215-248
\vskip 0.30cm
\item{[H]} H. Hasse, {\it Mathematische Abhandlungen}, Walter
de Gruyter, Berlin-New York, 1975.
\vskip 0.30cm
\item{[Mu]} D. Mumford, {\it Geometric Invariant Theory}, Springer-Verlag, 
Berlin
(1965)
\vskip 0.30cm
\item{[W]} A. Weil, {\it Sur les courbes alg\'ebriques et les vari\'et\'es qui
s'en d\'eduisent}, Herman, Paris (1948)
\vskip 0.30cm
\item{[We1]} L. Weng, New Non-Abelian Zeta Functions for Curves over Finite
Fields, Nagoya University, 2000
\vskip 0.30cm
\item{[We2]} L. Weng,  Riemann-Roch Theorem, Stability, and New Zeta
Functions for Number Fields, preprint, Nagoya, 2000
\eject
\vskip 2.0cm
\centerline {{\li Appendix:}\qquad {\we  Weierstrass  Groups}}
\vskip 0.5cm
\noindent
{\li 1. Weierstrass Divisors}

\noindent
(1.1) Let $M$ be a compact Riemann surface of
genus $g\geq 2$. Denote its degree $d$ Picard
variety by ${\rm Pic}^d(M)$. Fix a Poincar\'e
line bundle ${\cal P}_d$ on $M\times
{\rm Pic}^d(M)$. (One checks easily that  our
constructions  do not depend on a
particular choice of  Poincar\'e line bundle.)
Let $\Theta$ be the theta divisor of
${\rm Pic}^{g-1}(M)$, i.e., the image of the
natural map $M^{g-1}\to {\rm Pic}^{g-1}(M)$
defined by $(P_1,\dots,P_{g-1})\mapsto
[{\cal O}_M(P_1+\dots+P_{g-1})]$. Here $[\cdot]$
denotes the class defined by $\cdot$. We will
view the theta divisor as a pair $({\cal
O}_{{\rm Pic}^{g-1}(M)}(\Theta),{\bf 1}_\Theta)$
with ${\bf 1}_\Theta$ the defining section of $\Theta$ via
the structure exact sequence $0\to
{\cal O}_{{\rm Pic}^{g-1}(M)}\to {\cal
O}_{{\rm Pic}^{g-1}(M)}(\Theta)$.

Denote by $p_i$ the $i$-th projection of $M\times
M$ to $M$, $i=1,2$. Then for any degree $d=g-1+n$
line bundle on
$M$, we get a line bundle $p_1^*L(-n\Delta)$
on $M\times M$ which has relative $p_2$-degree
$g-1$. Here, $\Delta$ denotes the diagonal
divisor on $M\times M$. Hence, we get a
classifying map
$\phi_{L}: M\to {\rm Pic}^{g-1}(M)$ which makes the
following diagram commute:
$$\matrix{ M\times M&\to &M\times
{\rm Pic}^{g-1}(M)\cr
p_2\downarrow&&\downarrow \pi\cr
M&\buildrel\phi_{L}\over\to&
{\rm Pic}^{g-1}(M).}$$
One checks that there are  canonical
isomorphisms
$$\lambda_\pi({\cal P}^{g-1})\simeq {\cal
O}_{{\rm Pic}^{g-1}(M)}(-\Theta)$$ and
$$\lambda_{p_2}(p_1^*L(-n\Delta))\simeq
\phi_{L}^*{\cal O}_{{\rm
Pic}^{g-1}(M)}(-\Theta).$$ Here, $\lambda_\pi$
(resp. $\lambda_{p_2}$) denotes the
Grothendieck-Mumford cohomology determinant with
respect to $\pi$ (resp. $p_2$). (See e.g.,
[L].)

Thus, $\phi_{L}^*{\bf 1}_{\Theta}$ gives a
canonical holomorphic section of the dual of the
line bundle
$\lambda_{p_2}(p_1^*L(-n\Delta))$, which in turn
gives an effective divisor $W_L(M)$ on $M$, the
so-called {\it Weierstrass divisor associated to}
$L$.
\vskip 0.30cm
\noindent
{\it Example.} With the same notation as
above, take $L=K_M^{\otimes m}$ with $K_M$ the
canonical  line bundle of $M$ and $m\in
{\bf Z}$. Then we get an effective
divisor
$W_{K_M^{\otimes m}}(M)$ on $M$, which will be
called the $m$-{\it th  Weierstrass divisor}
associated to $M$. For simplicity, denote
$W_{K_M^{\otimes m}}(M)$ (resp.
$\phi_{K_M^{\otimes m}}$)
   by
$W_m(M)$ (resp. $\phi_m$).
\vskip 0.30cm
One checks easily that the degree of $W_m(M)$
is $g(g-1)^2(2m-1)^2$ and we have an
isomorphism ${\cal O}_M(W_m(M))\simeq
K_M^{\otimes g(g-1)(2m-1)^2/2}$.
Thus, in particular,
$$f_{m,n}:={{(\phi_m^*{\bf 1}_\Theta)^{\otimes
(2n-1)^2}}\over {(\phi_{n\ }^*{\bf
1}_\Theta)^{\otimes (2m-1)^2}}}$$ gives a
canonical meromorphic function on $M$ for all
  $m,n\in {\bf Z}$.
\vskip 0.30cm
\noindent
{\it Remark.} We may also assume that $m\in
{1\over 2}{\bf Z}$. Furthermore, this
construction has a relative
version as well, for which we assume that
$f:{\cal X}\to B$ is a semi-stable family of
curves of genus $g\geq 2$. In that case, we get
an effective divisor $({\cal O}_{\cal
X}(W_m(f)),{\bf 1}_{W_m(f)})$ and canonical
isomorphism
$$\eqalign{({\cal O}_{\cal
X}(W_m(f)),&{\bf 1}_{W_m(f)})\cr
\simeq&
({\cal O}_{\cal
X}(W_1(f)),{\bf 1}_{W_1(f)})^{\otimes
(2m-1)^2}\otimes ({\cal O}_{\cal
X}(W_{1\over 2}(f)),{\bf 1}_{W_{1\over
2}(f)})^{\otimes 4m(1-m)}.\cr}$$ The proof may be
given by using Deligne-Riemann-Roch theorem, which
in general, implies that we have the following
canonical isomorphism:
$$({\cal O}_{\cal X}(W_L(f)), {\bf 1}_{W_L(f)})\otimes
f^*\lambda_f(L)\simeq
L^{\otimes n}\otimes K_f^{\otimes n(n-1)/2}.$$
(See e.g. [Bur].) To allow   $m$ be a half
integer, we then should assume
that $f$ has a spin structure. Certainly,
without using spin structure,  a  modified
canonical isomorphism, valid  for integers, can
be given.
\vskip 0.30cm
\noindent
{\li 2. K-Groups}

\noindent
(2.1) Let $M$ be a compact
Riemann surface of genus $g\geq 2$. Then by the
localization theorem, we get the following exact
sequence for
$K$-groups
$$K_2(M)\buildrel\lambda\over\to K_2({\bf C}(M))\buildrel \coprod_{p\in 
M}\partial_p\over\to
\coprod_{p\in M}{\bf C}_p^*.$$

Note that the middle term may also be written as
  $K_2({\bf C}(M\backslash S))$ for any finite
subset $S$ of $M$, we see that naturally by the
theorem of Matsumoto,
the Steinberg symbol
$\{f_{m,n}, f_{m',n'}\}$ gives
a well-defined element in $K_2({\bf C}(M))$.
Denote the subgroup generated by all $\{f_{m,n},
f_{m',n'}\}$ with $m,n,m',n'\in {\bf Z}_{>0}$ in
$K_2({\bf C}(M))$ as $\Sigma(M)$.
\vskip 0.30cm
\noindent
{\bf Definition.} With the same notation as above,
the {\it first Weierstrass group} $W_I(M)$ of
$M$ is defined to be the $\lambda$-pull-back of
$\Sigma(M)$, i.e., the subgroup
$\lambda^{-1}(\Sigma(M))$ of
$K_2(M)$.
\vskip 0.45cm
\noindent
(2.2) For simplicity, now let $C$ be a regular
projective irreducible curve of genus $g\geq 2$
defined over
${\bf Q}$. Assume that $C$ has a semi-stable
regular module
$X$ over ${\bf Z}$ as well.
Then we have a natural
morphism $K_2(X)\buildrel\phi\over\to
K_2(M)$. Here $M:=C({\bf C})$.
\vskip 0.30cm
\noindent
{\bf Conjecture I.} {\it With the same notation as
above,  $\phi\big(K_2(X)\big)_{\bf Q}
=W_I(M)_{\bf Q}.$}
\vskip 0.45cm
\noindent
{\li 3. Generalized Jacobians}

\noindent
(3.1) Let $C$ be a projective, regular,
irreducible curve. Then for any effective divisor
$D$, one may canonically construct the so-called
generalized Jacobian $J_D(C)$ together with a
rational map
$f_D:C\to J_D(C)$.

More precisely, let $C_D$ be the group of classes
of divisors prime to $D$ modulo these which can
be written as ${\rm div}(f)$. Let $C_D^0$ be the
subgroup of $C_D$ which consists of all
elements of degree zero. For each $p_i$ in the
support of
$D$, the invertible elements modulo those
congruent to 1 (mod $D$) form an algebraic group
$R_{D,p_i}$ of dimension $n_i$, where $n_i$ is the
multiplicity of $p_i$ in $D$. Let $R_D$ be the
product of these $R_{D,p_i}$. One checks easily
that ${\bf G}_m$, the multiplicative group of
constants naturally embeds into $R_D$. It is a
classical result that we then have the short
exact sequence
$$0\to R_D/{\bf G}_m\to C_D^0\to J\to 0$$ where
$J$ denotes the standard Jacobian of $C$. (See
e.g., [S].) Denote $R_D/{\bf G}_m$ simply by
${\bf R}_D$.

Now the map $f_D$ extends naturally to a
bijection from $C_D^0$ to $J_D$. In this way
the commutative algebraic group $J_D$ becomes an
extension as algebraic groups of the standard
Jacobian by the group ${\bf R}_D$.
\vskip 0.30cm
\noindent
{\it Example.} Take the field of constants as
${\bf C}$ and $D=W_m(M)$, the $m$-th Weierstrass
divisor of a compact Riemann surface $M$ of
genus $g\geq 2$. By (1.1), $W_m(M)$ is
effective. So we get the associated generalized
Jacobian $J_{W_m(M)}$. Denote it  by
$WJ_m(M)$ and call it the $m$-th
{\it Weierstras-Jacobian} of $M$. For example, if
$m=0$, then $WJ_0(M)=J(M)$ is the standard
Jacobian of $M$. Moreover, one knows that
  the dimension of $R_{W_m(M),p}$ is at most
$g(g+1)/2$. For later use denote ${\bf
R}_{W_m(M)}$ simply by ${\bf R}_m$.
\vskip 0.45cm
\noindent
(3.2) The above construction works on any base
field as well. We leave the detail to the reader
while point out that if the curve is defined over
a field $F$, then its associated $m$-th
Weierstrass divisor is rational over the same
field as well. (Obviously, this is not true for
the so-called Weierstrass points, which behavior
in a rather random way.) As a consequence, by the
construction of the generalized Jacobian, we see
that the
$m$-th Weierstrass-Jacobians are also defined
over $F$. (See e.g., [S].)
\vskip 0.45cm
\noindent
{\li 4. Galois Cohomology Groups}

\noindent
(4.1) Let $K$ be a perfect field, $\overline K$
be an algebraic closure of $K$ and $G_{\overline
K/K}$ be the Galois group of $\overline K$ over
$K$. Then for any $G_{\overline K/K}$-module $M$,
we have the Galois cohomology groups
$H^0(G_{\overline K/K},M)$ and $H^1(G_{\overline
K/K},M)$ such that if $$0\to M_1\to M_2\to M_3\to
0$$ is an exact sequence of $G_{\overline
K/K}$-modules, then we get  a natural long exact
sequence
$$\eqalign{0\to &H^0(G_{\overline
K/K},M_1)\to
H^0(G_{\overline
K/K},M_2)\to H^0(G_{\overline
K/K},M_3)\cr
&\to H^1(G_{\overline
K/K},M_1)\to H^1(G_{\overline
K/K},M_2)\to H^1(G_{\overline K/K},M_3).\cr
}$$
Moreover, if $G$ is a subgroup
of $G_{\overline K/K}$ of finite index or a
finite subgroup, then
$M$ is naturally a $G$-module. This leads a
restriction map on cohomology
${\rm res}: H^1(G_{\overline K/K},M)\to
H^1(G,M)$.
\vskip 0.30cm
\noindent
(4.2) Now let $C$ be a projective,
regular irreducible curve defined over a number
field $K$. Then for each place $p$ of $K$, fix an
extension of $p$ to $\overline K$, which then
gives
  an embedding
$\overline K\subset \overline K_p$ for the
$p$-adic completion $K_p$ of $K$ and a
decomposition group $G_p\subset G_{\overline
K/K}$.

Now apply the  construction in (3.1) to the short
exact sequence
$$0\to {\bf R}_m\to WJ_m(C)\to J(C)\to 0$$ over
$K$. Then we have the following
long exact sequence
$$\eqalign{0\to& {\bf R}_m(K)\to WJ_m(K)\to
J(K)\cr
&\to H^1(G_{\overline K/K}, {\bf R}_m(K))\to
H^1(G_{\overline K/K}, WJ_m(K))\buildrel
\psi\over\to H^1(G_{\overline K/K},J(K)).\cr}$$

Similarly, for each place $p$ of $K$, we have  the
following exact sequence
$$\eqalign{0\to& {\bf R}_m(K_p)\to WJ_m(K_p)\to
J(K_p)\cr
&\to H^1(G_p, {\bf R}_m(K_p))\to
H^1(G_p, WJ_m(K_p))\buildrel
\psi_p\over\to
H^1(G_p,J(K_p)).\cr}$$

Now the natural inclusion $G_p\subset
G_{\overline K/K}$ and $\overline
K\subset\overline K_p$ give restriction maps on
cohomology, so we  arrive at a natural
morphism
$$\Phi_m:\psi\Big(H^1(G_{\overline K/K},
WJ_m(K))\Big)\to
\prod_{p\in M_K}\psi_p\Big(H^1(G_p,
WJ_m(K_p))\Big).$$
Here $M_K$ denotes the set of all places over $K$.
\vskip 0.30cm
\noindent
{\bf Definition.} With the same notation as above,
{\it the second Weierstrass group} $W_{II}(C)$ of
$C$ is defined to be the subgroup of
$H^1(G_{\overline K/K}, J(C)(K))$
generated by all ${\rm Ker}\,\Phi_m$, the kernel
of $\Phi_m$, i.e.,
$W_{II}(C):=\langle {\rm Ker}\,\Phi_m:m\in {\bf
Z}_{>0}\rangle_{\bf Z}.$
\vskip 0.30cm
\noindent
{\bf Conjecture II.} {\it With the same notation
as above, the second Weierstrass group $W_{II}(C)$
is   finite.}
\vskip 0.45cm
\noindent
{\li 5.  Deligne-Beilinson Cohomology}

\noindent
(5.1) Let $C$ be a projective regular curve of
genus $g$. Let $P$ be a finite set of
$C$. For simplicity, assume that all of them
are defined over {\bf R}. Then
we have the associated Deligne-Belinsion
cohomology group
$H_{\cal D}^1(C\backslash P,{\bf R}(1))$ which
leads to the following short exact sequence:
$$0\to {\bf R}\to H_{\cal D}^1(C\backslash P,{\bf
R}(1))\buildrel{\rm div}\over\to {\bf R}[P]^0\to
0$$ where ${\bf R}[P]^0$ denotes  $({\rm the\
group\ of\ degree\ zero\ divisors\ with\ support\
on\ } P)_{\bf R}$.

The standard cup product  on
Deligne-Beilinson cohomology leads to a
well-defined map:
$$\cup:H_{\cal D}^1(C\backslash P,{\bf
R}(1))\times H_{\cal D}^1(C\backslash P,{\bf
R}(1))\to H_{\cal D}^2(C\backslash P,{\bf
R}(2)).$$

Furthermore, by Hodge theory, there is a canonical
short exact sequence
$$0\to H^1(C\backslash P,{\bf R}(1))\cap
F^1(C\backslash P)\to H_{\cal D}^2(C\backslash
P,{\bf R}(2))\buildrel {p_{\cal D}}\over \to
  H_{\cal D}^2(C,{\bf R}(2))\to 0$$ where
$F^1$ denotes the $F^1$-term of the Hodge
filtration on $H^1(C\backslash P,{\bf C})$.

All this then leads to a well-defined morphism
$$[\cdot,\cdot]_{\cal D}:\wedge^2 {\bf R}[P]^0\to
H_{\cal D}^2(C,{\bf R}(2))=H^1(C,{\bf R}(1))$$
which make the associated diagram coming from the
above two short exact sequences commute. (See e.g.
[Bei].)

\noindent
(5.2) Now applying the above construction with
$P$ being the union of the  supports of $W_1$,
$W_m$ and $W_n$ for $m,n>0$. Thus for fixed
$m,\,n$, in ${\bf R}[P]^0$,  we get two
elements
${\rm div}(f_{1,m})$ and ${\rm div}(f_{1,n})$.
This then gives $[{\rm div}(f_{1,m}), {\rm
div}(f_{1,n})]_{\cal D}\in H^1(X,{\bf R}(1)).$

\noindent
{\bf Lemma.} {\it For any holomorphic
differential 1-form $\omega$ on $C$, we have
$$\eqalign{~&\langle [{\rm div}(f_{1,m}), {\rm
div}(f_{1,n})]_{\cal D},\omega\rangle
:=-{1\over{2\pi {\sqrt -1}}}\int [{\rm
div}(f_{1,m}), {\rm div}(f_{1,n})]_{\cal
D}\wedge\bar\omega\cr =&-{1\over{2\pi {\sqrt -1}}}
\int g({\rm div}(f_{1,m}),z)d g({\rm
div}(f_{1,n}),z)\wedge\bar\omega,\cr}$$ Here
$g(D,z)$ denotes the Green's function of $D$ with
respect to any fixed normalized (possibly
singular) volume form of quasi-hyperbolic type.}

\noindent
{\bf Proof.} A simple argument by using the Stokes
formula.

\noindent
(4.3) With exactly the same notation as in (4.2), then in
$H^1(X,{\bf R}(1))$ we get a collection of
elements $[{\rm div}(f_{1,m}), {\rm
div}(f_{1,n})]_{\cal D}$ for $m,n\in {\bf
Z}_{>0}$.

\noindent
{\bf Definition.}  With the same notation as
above, assume that $C$ is defined over {\bf Z}.
Define the {\it --first quasi-Weierstrass group}
$W_{-I}'(C)$ of $C$ to be the subgroup of
$H^1(X,{\bf R}(1))$ generated by $[{\rm div}(f_{1,m}), {\rm
div}(f_{1,n})]_{\cal D}$ for all $m,n\in {\bf
Z}_{>0}$ and call  $W_{-I}'(C)_{\bf Q}$
the {\it --first Weierstrass group} $W_{-I}(C)$ of
$C$. That is to say,
$W_{-I}(C):=\langle [{\rm div}(f_{1,m}), {\rm
div}(f_{1,n})]_{\cal D}: m,n\in {\bf
Z}_{>0}\rangle_{\bf Q}.$
\vskip 0.30cm
\noindent
{\bf Conjecture III.} {\it With the same notation
as above, $W_{-1}(C)_{\bf R}$ is the full space,
i.e. equals to
$H^1(X,{\bf R}(1))$.}
\vskip 0.30cm
\noindent
That is to say, we believe Weierstrass divisors
will give a new rational structure for $H^1(X,{\bf
R}(1))$. Furthermore, we believe that the
corresponding regulator will give the leading
coefficient of the $L$-function of $C$ at $s=0$,
up to rationals.
\vskip 0.45cm
\noindent
\centerline {\bf REFERENCES}
\vskip 0.30cm
\item{[Bei]} A. A. Beilinson, Higher
regulators and values of $L$-functions, J. Soviet
Math. 1985
\vskip 0.30cm
\item{[Blo]} S. Bloch, {\it Higher regulators,
algebraic $K$-theory and zeta-function of
elliptic curves}, UC Irvine, 1978
\vskip 0.30cm
\item{[Bur]} J.-F. Burnol, Weierstrass points
on arithmetic surfaces, Invent. Math. 1992,
\vskip 0.30cm
\item{[L]} S. Lang, {\it Introduction to
Arakelov theory}, Springer-Verlag, 1988
\vskip 0.30cm
\item{[S]} J. P. Serre,  {\it Groupes
algebriques et corps de classes}, Hermann, 1975
\end